\documentclass{article}

\usepackage{amscd, amssymb}

\usepackage{graphicx}
\usepackage{psfrag}

\newtheorem{thm}{Theorem}

\newtheorem{cor}{Corollary}

\newcommand{\ignore}[1]{\index{ignore}}

\renewcommand{\O}{{\mathcal O}}

\renewcommand{\d}{\delta}
\newcommand{\D}{\Delta}
\renewcommand{\O}{\Omega}
\newcommand{\w}{\omega}

\newcommand{\E}{{\mathbb E}}
\newcommand{\isom}{\cong}
\newcommand{\text}{\mbox}

\begin{document}
\title{The number of possibilities for random dating\footnote{This note appeared in \emph{Journal of Combinatorial Theory, Series A}, vol.~115 (2008), pp.~1265--1271.}}

\author{
    {Aaron Abrams}\footnote{
    Department of Mathematics and Computer Science, Emory University, Atlanta GA 30322, USA,
    \tt{abrams@mathcs.emory.edu}}
    \ and {Rod Canfield}\footnote{
    Department of Computer Science, University of Georgia, Athens GA 30602, USA,
    \tt{erc@cs.uga.edu}}
    \ and {Andrew Granville}\footnote{
    D\'epartement de Math\'ematiques et Statistique, Universit\'e de Montr\'eal, CP 
    6128 succ Centre-Ville, Montr\'eal QC  H3C 3J7, Canada,
    \tt{andrew@dms.umontreal.ca}} \\
}
\date{}
\maketitle

\begin{abstract}
Let $G$ be a regular graph and $H$ a subgraph on the same vertex set.
We give surprisingly compact formulas for the number of copies of $H$
one expects to find in a random subgraph of $G$.
\end{abstract}

\section{Introduction}

There are $n$ boys and $n$ girls who apply to a computer dating
service, which randomly picks a boy and a girl and then introduces
them. It does this again and again until everyone of the boys and
girls has been introduced to at least one other person. The
service then organizes a special evening at which everyone dates
someone to whom they have been previously introduced. In how many
different ways can all of the boys and girls be matched up?

Translating this question into the language of graph theory, we
select edges at random from the complete bipartite graph $K_{n,n}$
until the subgraph created by these edges, $G$, has minimum degree
1. We then ask how many perfect matchings (subgraphs on the same
vertex set where every vertex has degree 1) are contained in $G$?

We know that, with probability $1-o(1)$, there will be at least
one possible matching [1] at the time that the last person is
finally introduced to someone, and that there will have been a
total of $\sim n\log n$ introductions made [2]. Now, note that if
the $j$th girl has been introduced to $b_j$ boys then the
number of possible matchings is certainly $\leq b_1b_2\dots
b_n \leq (k/n)^n$ where $k=b_1+b_2+\dots +b_n$, by the
arithmetic-geometric mean inequality. In particular we expect no
more than $(\{ 1+o(1)\} \log n)^n$ ways of matching up the boys
and girls, by the time   the last person is finally introduced to
someone. It thus came as a surprise to us, when we did the
calculation, to discover that one expects there to be far more,
in fact something like $(n/4\text{e})^n$ ways of matching up the 
boys and girls!

The fault in the above reasoning comes in believing that what we
expect to happen in the most likely case (in which there have been
$\sim n\log n$ introductions made by the time the last person is
finally introduced to someone) yields what we expect to happen in
general. In fact the largest contribution to the expectation takes
place in the very rare situation that the computer dating service
somehow neglects to make any introductions involving one sad
participant  until about $n^2/2$ introductions have already been made
involving only the others.

To see this, first note that the probability that the $n$th boy
has not been introduced to anyone by the time $k=[n^2/2]$ random
introductions have been made is\footnote{For this and later
estimates, note that $\log(A(A-1)\dots
(A-(n-1))/A^n)=-n^2/2A +O(1/n)$ for $A\gg n^2$.}
$$
 \frac{n^2-n}{n^2} \cdot \frac{n^2-n-1}{n^2-1}  \cdot
\frac{n^2-n-2}{n^2-2} \cdots    \frac{n^2-n-k+1}{n^2-k+1}\sim
\frac{\text{e}^{-1/2}}{2^n} .
$$
With probability $1+o_n(1)$, each person other than the $n$th boy has, by
now, been introduced to $\sim n/2$ people (here  $o_n(1)$ represents a function that $\to 0$ as $n\to \infty$).  
Suppose that the
$(k+1)$st introduction involves the $n$th boy and girl. For any of
the $(n-1)!$ possible matchings involving all of the boys and
girls in which the $n$th boy and girl are matched up, the
probability that each of the pairings in that matching has
already been introduced is close to $1/2$, so we might guess that
the probability that the matching can occur with the introductions
already made is around $1/2^{n-1}$. Thus we might expect that the
contribution to the expectation when $k=[n^2/2]$ is around 
$(n-1)!\times (\text{e}^{-1/2}/2^n) \times (1/2^{n-1}) \approx
(n/4\text{e})^n$.

The main point of this article is to make this discussion precise
with complete proofs, and in some generality. The exact
formulas obtained are surprisingly compact.

\section{Statement of results}

Let $G_{\w}$ be a random subgraph of $K_{n,n}$ formed by randomly
adding edges until every  vertex has degree at least one. Various
things are known about such $G_{\w}$; for instance $G_{\w}$ is
almost surely connected  and $G_{\w}$ almost surely contains a
matching [1]. We ask  how many matchings does $G_{\w}$ contain?

One can ask the same question for random subgraphs $G_{\w}$ of
$K_{2n}$, or indeed of any other graph $G$.  Likewise instead of
matchings one could count the expected number of occurrences of
any other prescribed subgraph $H$ on the same vertex set as $G$.
Of course, before any copy of $H$ can appear as a subgraph of
$G_{\w}$, it is necessary that the minimum degree of $G_{\w}$ be
at least as great as the minimum degree of $H$. Note also that $G$
itself may not contain any subgraphs isomorphic to $H$; thus it
makes sense to investigate the proportion of subgraphs of $G$
which are isomorphic to $H$ that actually occur as a subgraph of
$G_\w$. Define
$$
\E(H\subseteq G):= \E
\left( \frac {\#\{J\subseteq G_{\w} \ : \ J\isom H \} } {
\#\{J\subseteq G : J\isom H\} } \right).
$$

Thus let $G$ and $H$ be graphs with $|V(G)|=|V(H)|$, and let $\d\geq 1$ be
the minimum degree of $H$.  Let $\O$ be the set of permutations
of the $m$ edges of $G$.  For $\w=(e_{\w(1)},\ldots,e_{\w(m)})\in\O$,
let $G_{\w}^{(i)}$ be the graph with vertex set $V(G)$ and edge set
$\{e_{\w(1)},\ldots,e_{\w(i)}\}$, and let
$$k(\w)=\min \{i: G_{\w}^{(i)} \text{ has minimum degree } \d\}.$$
Then define $G_{\w}=G_{\w}^{(k(\w))}$.  The question now becomes:  for
random $\w\in\O$, what fraction of those subgraphs of $G$ which are
isomorphic to $H$ are contained in $G_{\w}$?

%What is the expected number of distinct subgraphs
%of $G_{\w}$ which are isomorphic to $H$?

Our main theorem answers this question in the case that $G$ is regular.

\begin{thm}
Suppose $G$ is $d$-regular, and let $H$ be a subgraph of $G$ with
the same vertex set, of minimal degree   $\d\geq 1$. Let $h$ be
the number of edges of $H$, and let $\D=d-\d$. Then
$$\E(H\subseteq G)=
\ \frac 2 {{{h+\D}\choose h}} - \frac 1 {{{h+2\D}\choose h}} 
.$$
\end{thm}

We can apply this to our original problem.

\begin{cor}
In a random subgraph of $K_{n,n}$ formed by adding edges until each vertex
has degree at least one, the expected number of complete matchings is
$$n!\left(\frac {2} {{{2n-1}\choose n}} - \frac {1} {{{3n-2}\choose n}}
\right).$$
\end{cor}

{\bf Proof.}
In the notation of the theorem, $d=n$, $\d=1$, $h=n$, and $\D=n-1$.
Apply the theorem, and note that there are $n!$ distinct complete
matchings in $K_{n,n}$.
\hfill $\square$ \\

If boys can date boys, and girls can date girls, we can still ask
the same questions, and again obtain what is to us a surprising answer.

\begin{cor}
Let $G$ be any $d$-regular graph on $2n$ vertices, and form
$G_{\w}$ by randomly choosing edges of $G$ until each vertex has
degree at least one. Of all the complete matchings in $G$, the
fraction expected to occur in $G_{\w}$ is
$$\frac {2}{{{n+d-1}\choose n}} - \frac {1} {{{n+2d-2}\choose n}}.$$
\end{cor}

\begin{cor}
Let $G$ be a $d$-regular graph on $n$ vertices, and form $G_{\w}$
by randomly choosing edges of $G$ until each vertex has degree at
least two. Of all the Hamiltonian cycles in $G$, the fraction
expected to occur in $G_{\w}$ is
$$\frac 2 {{{n+d-2}\choose n}} - \frac 1 {{{n+2d-4}\choose n}}.$$
\end{cor}

The proof of the main theorem uses a simple inclusion-exclusion
argument. It is not hard to prove a result for more general graphs
$G$ and $H$ as indicated in the proof, though this would be
complicated to state. In [3] McDiarmid proves these results in the
special case that $G=K_n$ though with more examples of desired
subgraphs (that is, other than just for matchings and for Hamiltonian
cycles).  In our proof, in section 3, we are able to take any $G$, and
any desired subgraph $H$, so long as it has as many vertices as $G$.
Our proof has similarities to that in [3], though it is not entirely clear
how to generalize [3] directly to recover the general results presented
in section 3.

There is a hypergraph version whose proof is
essentially the same but whose statement is correspondingly more
complicated. Thus we will state two special cases.

First, let $G=K_{n,n,\dots,n}$ be a complete $r$-partite graph with
vertex set $V(G)$.  Let $E$ be the set of $r$-element subsets of
$V(G)$ which span $K_r$ subgraphs of $G$ (so $|E|=n^r$).  Form an
$r$-uniform hypergraph (on the same vertex set as $G$) by choosing
random elements of $E$ as edges  until each vertex is contained in
at least one edge; i.e. choose a permutation $\w$ of $E$ and let
$G_{\w}$ be defined as it was before.

A \emph{matching} in such a hypergraph $G_{\w}$ is a set $M$ of edges such
that each vertex of $G$ is in exactly one edge of $M$.  (Equivalently
a matching is a copy, inside $G_{\w}$, of the hypergraph $H$ consisting of
$n$ disjoint edges, each of size $r$.)

\begin{thm}
With $G=K_{n,n,\dots,n}$ the complete $r$-partite graph, and with $G_{\w}$
as above, the expected number of matchings in $G_{\w}$ is
$$(n!)^{r-1} \sum_{i=1}^r (-1)^{i-1} \frac {{r\choose i}}
{{{n^r - (n-1)^i n^{r-i} + n - 1}\choose n}} .$$
\end{thm}

The factor $(n!)^{r-1}$ is  the total number of matchings in $G$.
The case $r=2$   agrees with Corollary 1. For $r\geq 3$ we have
that $n^r - (n-1)^i n^{r-i} + n - 1=in^{r-1}(1+O(1/n))$ so that
${{n^r - (n-1)^i n^{r-i} + n - 1}\choose n} \asymp
(in^{r-1})^n/n!$. Thus the $i=1$ term is the main term in the sum,
and so the
sum can be estimated as $\{ rc_r +O_r(1/n)\} (n!)^r/n^{(r-1)n}$
where $c_3=\text{e}^{-1/2}$ and $c_r=1$ for $r\geq 4$. A more
crude estimate is $\asymp_r n^{O_r(1)} (n/\text{e}^r)^n$ 
(where ``$O_r(1)$'' is in place of a function that is bounded in terms of $r$ only, and where ``$A  \asymp_r B $'' means that $A/B$ is bounded above and below by positive constants that depend only on $r$).

For our second hypergraph theorem, let $G$ be the complete graph
$K_{rn}$. Let $E$ be the set of all $r$-element subsets of $V(G)$.
We form an $r$-uniform hypergraph $G_{\w}$ by taking edges from
$E$ in some order $\w$, stopping when each vertex has degree at
least one.  Again, we find the expected number of matchings in
such a hypergraph $G_{\w}$.

\begin{thm}
With $G=K_{rn}$ and with $G_{\w}$ as above, the expected number of matchings
in $G_{\w}$ is
$$\frac {(rn)!}{(r!)^n n!} \sum\limits_{i=1}^r
(-1)^{i-1} \frac {{r\choose i}}{{
{ {{nr} \choose r} - {{nr-i}\choose r} +n -1 }\choose n }}.$$
\end{thm}

Again, the factor preceding the sum is the total number of
matchings in $G$. The case $r=2$   agrees with Corollary 2 with
$G=K_{2n}$. For $r\geq 3$ we have that if $i\geq 2$ then  ${{nr}
\choose r} - {{nr-i}\choose r} \geq ({{nr} \choose r} -
{{nr-1}\choose r}) (2 - \frac{(r-1)}{(nr-1)})$, and thus the
$i\geq 2$ terms have magnitude $\ll 1/2^n$ times the magnitude of
the $i=1$ term. Therefore the sum can be estimated as $\{ rc'_r
+O_r(1/n)\}   ( ((n-1)r)! n)^n/(rn)!^{n-1}$ where
$c'_3=\text{e}^{-1/9}$ and $c'_r=1$ for $r\geq 4$. A more crude
estimate is  $\asymp_r n^{O_r(1)} (n/\text{e}^r)^n$.

The proofs of theorems 2 and 3 are the obvious
generalizations of the proof of theorem 1. Thus we have restricted
ourselves to sketching the proof of theorem 2, and leaving the
proof of theorem 3 to the enthusiastic reader.

\section{Proof of theorem 1}

We are given graphs $G$ and $H$ with $V(G)=V(H)$ and where
$\min_{v\in V(G)} \deg_G(v) \geq \min_{v\in V(G)} \deg_H(v)$.  For 
now we need not assume that $G$ is regular.  We use the notation
$\d$, $\O$, $k(\w)$, $G_{\w}$, and $h$ as defined in the previous section.

By definition
$$\E(H\subseteq G)= \frac 1{\#\{J\subseteq G : J\isom H\}} \
\sum_{{J\subseteq G} \atop {J\isom H}}
\Pr  (e_1,\ldots,e_{k(\w)} \supseteq J:\ \w\in\O).$$

Now we fix $J\subseteq G$ with $J\cong H$ and evaluate the probability that $J\subseteq e_1,\ldots,e_{k(\w)}$ as we vary over $\w\in\O$. 
Fix $\w$, and let
$u$ and $v$ be the endpoints of the edge $e_{k(\w)}$.  By definition
of $k(\w)$, either $\deg_{G_{\w}}(u)=\d$ or $\deg_{G_{\w}}(v)=\d$ (or
both).

Now, $J\subseteq G_{\w}$ with $e_{k(\w)}=uv$ and
$\deg_{G_{\w}}(u)=\d$ if and only if certain edges of $G$ appear
in the correct order in $\w$:
\begin{itemize}
\item
first, the edges of $J$ other than $uv$ (in some order);
\item
then $uv$;
\item
then the edges of $G$ which contain $u$ other than those in $J$.
\end{itemize}
Note that the location in $\w$ of the other edges of $G$ is irrelevant.

The probability that the edges appear in this order is
$$\frac{ (h-1)! 1! (\deg_G (u) - \d)! }{ (h + \deg_G (u) - \d)! } \ ,$$
where $h$ is the number of edges of $H$.

A similar count holds if $\deg_{G_{\w}}(v)=\d$.  However, we have now
double counted the cases in which $\deg_{G_{\w}}(u) = \d = \deg_{G_{\w}}(v).$
By similar reasoning, the probability of this is
$$\frac{ (h-1)!1!(\deg_G(u) + \deg_G(v) - 2\d)! }{ (h+\deg_G(u)+\deg_G(v)-2\d)!}
.$$
Thus, given that the last edge $e_{k(\w)}$ is the edge $uv$, we have
$$\Pr (e_1,\ldots,e_{k(\w)}\supseteq J, \ \text{ and }
e_{k(\w)}=uv:\ \w\in\O)= \frac 1 h \left( \frac 2 {{{h+\D(u)}\choose h}} -
        \frac 1 {{{h+\D(u)+\D(v)} \choose h}} \right),$$
where $\D(u)=\deg_G(u)-\d$. Note that the only dependence of the right side on $J$ is the choice of the edge $uv$ from $J$.
We deduce immediately that 
$$
\Pr_{\w\in\O}(e_1,\ldots,e_{k(\w)}\supseteq J:\ \w\in\O)= 
\sum_{uv \in J} \ 
\frac 1 h \left( \frac 2 {{{h+\D(u)}\choose h}} -
        \frac 1 {{{h+\D(u)+\D(v)} \choose h}} \right).$$
Now if $G$ is $d$-regular then $\D(u)=d-\d$ for all vertices $u$,
and our sum is over $h$ edges, so that
$$\Pr_{\w\in\O}(e_1,\ldots,e_{k(\w)}\supseteq J)=
\frac 2 {{{h+\D}\choose h}} - \frac 1 {{{h+2\D} \choose h}}$$
unconditionally.
\hfill $\square$ \\

\bigskip

\noindent{\em Remark}:  The proof above yields, when $G$ is not 
regular, that the expected number of copies of $H$ in $G_{\w}$ is
$$
\frac 1 h \sum_{uv \text{ is an edge of } G}
\left( \frac 2 {{{h+\D(u)}\choose h}} -
        \frac 1 {{{h+\D(u)+\D(v)} \choose h}} \right) 
        \# \{ J\subseteq G:\ J\cong H,\ \text{and} \ uv\in J\} .
$$

\section{Proof of theorem 2}

Denote the vertices of $K_{n,\ldots,n}$ by $v_{i,j}$ where $1\leq i\leq r$
and $1\leq j \leq n$.  Each edge of $E$ has one vertex $v_{i,j(i)}$ for
each $i$; we abbreviate the edge as $(j(1),\ldots,j(r))$.

To prove theorem 2 we apply the same technique as in the proof of theorem 1.
The expected number of matchings in $G_{\w}=\{e_1,\ldots,e_{k(\w)}\}$ is
the total number of possible matchings (which is $(n!)^{r-1}$) times the
probability that any given matching occurs in $G_{\w}$.  We may compute
the latter probability for the ``diagonal" matching
$M=\{(1,\ldots,1),\ldots,(n,\ldots,n)\}$ as follows:

\begin{eqnarray*}
\Pr (G_{\w}\supseteq M) &=& n \Pr \left(G_{\w}\supseteq M\ \text{ and }
e_{k(\w)}=(n,\ldots,n) \right) \\
&=& n \sum\limits_{i=1}^r (-1)^{i-1} {r\choose i} \Pr
(G_{\w}\supseteq M,\  e_{k(\w)}=(n,\ldots,n),\ \\
& &    \hskip1in   \text{and }
v_{1,n},\ldots,v_{i,n}\ \text{ all have degree one in}\ G_{\w}), \\
\end{eqnarray*}
where the last equality follows by inclusion-exclusion, and the
symmetry of $G$.  We can determine the latter probability as
before, namely, the event occurs if and only if the relevant edges
occur in $\w$ in the correct order:
\begin{itemize}
\item
first, the edges $\{(i,\ldots,i)\ :\ 1\leq i\leq n-1 \}$ (in some order);
\item
then $(n,\ldots,n)$;
\item
then any edge (other than $(n,\dots,n)$) containing $v_{j,n}$ for some $j$,
$1\leq j \leq i$.
\end{itemize}

Again, the other edges can occur anywhere in $\w$.  Thus the desired
probability is
$$\frac {(n-1)! (n^r-(n-1)^i n^{r-i}-1) !} { (n^r-(n-1)^i n^{r-i} +n-1)!}
= \frac 1 n { {n^r-(n-1)^i n^{r-i} +n-1} \choose n}^{-1}.$$

This establishes theorem 2.

\section{Contributions to the expectation in Theorem 1}

To better understand the result in Theorem 1, we now determine the
contribution to the expectation of the sum, over all $\w$ with
$k(\w)=k$, of the number of copies of $H$ contained in
$e_1,\ldots,e_k$. Let $E$ be the total number of edges in $G$.
Just as in the proof of Theorem 1 the main term is given by $2h$
times the number of sets $e_1,\ldots,e_k$ which contain $J$, with
$e_k=uv$, where the other edges of $G$ which contain $u$ lie
outside $e_1,\ldots,e_k$. Once these edges are chosen, the number
of possibilities for the $k-h$ edges in $\{e_1,\ldots,e_k\}
\setminus J$ is ${E-h-\Delta \choose k-h}$. The number of possible
orderings of the edges in $G$ given the above restrictions is
$(k-1)!1!(E-k)!$. By an analogous calculation when $u$ and $v$
both have degree $\delta$ in $G_\w$, we have that the sum, over
all $\w$ with $k(\w)=k$, of the number of copies of $H$ contained
in $e_1,\ldots,e_k$, is
$$
h (k-1)!1!(E-k)!  \left\{ 2 {E-h-\Delta \choose k-h} -
{E-h-2\Delta \choose k-h} \right\}.
$$
If $h=o(k)$ and $k=o(E)$ then this is $E! (k/E)^h
\text{exp}(O(k^2/E+h^2/k))$. It is maximized when $k$ is a little larger, 
namely when $k \sim
hE/(\Delta+h)$ (and note that there are no more than $E-h+1$
possible values for $k$).  Since $\Delta\leq d-1 \leq
|V(G)|-2=|V(H)|-2< \delta |V(H)|\leq 2h$, we have 
$E\geq hE/(\Delta+h)\geq E/3$.

In the case considered in Corollary 1, $h=n,\Delta=n-1$ and
$E=n^2$, so the main contribution to the expectation occurs when
$k\sim n^2/2$, far more edges than when $k\sim n\log n$ which is
the value we expect for $k(\w)$.  Likewise, for matchings in $K_{2n}$
the biggest contribution occurs with $k=2n^2/3$ rather than
the expected $k\sim n\log n$, and for Hamiltonian cycles in $K_n$
the biggest contribution comes from $k=n^2/4$, not 
$k\sim (n/2)\log n$.


\newpage

\begin{thebibliography}{[K-K-M]}


\bibitem[1]{boll} B\'ela Bollob\'as, {\em Random Graphs}, Cambridge
Studies in Advanced Mathematics {\bf 73} (2001).


\bibitem[2]{et} Paul Erd\H os and Alfr\'ed R\'enyi, {\em On random graphs I},
Publ. Math. Debrecen  {\bf 6} (1959) pp. 290--297.


\bibitem[3]{mcd} Colin McDiarmid, {\em Expected Numbers at Hitting Times},
J. Graph Theory {\bf 15} (1991) pp. 637--648.

\end{thebibliography}
\end{document}